\documentclass{notices}
\usepackage{amsfonts,amssymb,amsmath,amscd}
\usepackage{xcolor}
\usepackage{graphicx}
\usepackage{subcaption}

\title{
Extreme cases: Math education within the U.S. prison system
}

\author{
  Katherine J. Pearce
  \affil{
    KJP is a Peter O'Donnell Jr. Postdoctoral Fellow in the Oden Institute at the University of Texas at Austin and volunteers with the Texas Prison Education Initiative. Her email address is katherine.pearce@austin.utexas.edu.
    }
}

\begin{document}

\maketitle

Pulling into a visitor parking space for the first time, I watch a stray dog hesitantly making her way toward the barbed wire fences. 
The surroundings are austere and unforgiving. Lockhart Correctional Facility is about an hour’s drive from the city of Austin and miles from the nearest town center, situated, intentionally, in the middle of nowhere. 
I wonder briefly how the dog ended up here and feel a pang of sadness that I have nothing to offer her. 
It will not be the last time I struggle with these feelings in this parking lot; they only intensify when I think about my intelligent and resourceful students being left here to fend for themselves. 

I pocket my ID and car keys, grab the box of printouts I made for class, and head toward the barbed wire gate. 
With sunlight still streaming through the lobby’s glass doors behind me, the first checkpoint feels similar to airport security. 
I empty my pockets, place everything on a table to be searched, remove my belt and shoes, and walk slowly through the metal detector. 
A female correctional officer (CO) outlines my body with her wand before patting me down, paying special attention to the bottoms of my feet, while another CO searches my class materials. 
He warns me that highlighters and papers with “too much ink” are not permitted, but thankfully all my materials are allowed through today. 
I gather them up and head toward the next security checkpoint, where any glimpse of daylight or sense of familiarity is gone. 

The main guard station remotely buzzes me through a door after surveilling me on camera, and I hand over my ID in exchange for a badge after signing in with them. 
I am buzzed through one last door and suddenly step out into prison, narrowly avoiding the lines of people moving seamlessly in and out of the cafeteria.
I hear my pulse in my ears until I notice that they are all smiling and welcoming me.
Someone exclaims, “God, I miss wearing jeans!” 
We all laugh.
Another woman approaches me to show me to my classroom, and we walk down a hallway under fluorescent lighting and a large hemispherical mirror that shows everything around the corner. 
“We are so grateful and excited you are here!” she tells me, and I feel the same. 

When we arrive at the classroom, she apologetically explains that I will have to wait while the COs release my students to come to class. 
Although class technically begins at 6 p.m., we usually don’t start until 6:30, but I enjoy the wait. 
That half hour is the most important part of teaching for me; without a phone or any electronics to distract me, I wind up studying all the art and handouts displayed on the walls.
This classroom is also used by vocational and re-entry programs, and the exercises that the students complete are somber reminders of the real stakes here. 
“My goal is to get sober for my children.” 
“I want to learn how to forgive myself and earn forgiveness from others.” 
“When I get out, I will start a new career with my certifications to support my family.” 
These messages reiterate to me that I am primarily here to build their confidence in themselves and their problem-solving abilities. This confidence gives them a sense of freedom even in incarceration, as a formerly incarcerated friend will later share with me.

The twelve students in my class trickle in, and we introduce ourselves. 
Even though I learned about average state prison populations during the Texas Prison Education Initiative's (TPEI) orientation, I am still surprised at the large variance in age, from about early 20s to mid 60s.
I feel a flutter of anxiety about the pace I've decided on for my lectures. 
The students need to pass my course as a prerequisite for any credit-bearing TPEI math course through UT Austin, so I want to focus on building their abstract problem-solving abilities.
In particular, I have decided to spend a lot of time developing their intuition in subject areas they had likely seen before, like integer addition; I am worried the students will be bored with my decision. 

``So, just out of curiosity, what do y'all think about math in general?'' I ask, quickly adding, ``No judgment here, I do it for a living and still feel love-hate about it.''
We laugh, but it's hard to keep smiling encouragingly when I hear about their previous experiences in math classes. 
The students are diverse in age, race, ethnicity, and sexual orientation, but they unfortunately had one major thing in common: all of them proclaimed to be ``bad at math,'' and many had dropped or failed out of school specifically because of math. \footnote{According to the Bureau of Justice Statistics, fewer than 4 out of 10 incarcerated people have completed high school, versus 9 out of 10 in the general population.
The average incarcerated person in state prison is 39 years old with a 10th grade education. \cite{PPI22}} 

I tell them my plan for the course: we are going to start over with math and look at it from a new perspective. 
Anything that they need to know for the course, I will teach them, and it would actually be preferable to wipe away any memory of subjects where they had unsuccessful learning outcomes.
Even things they have seen before, like integer addition and multiplication, I want us to consider with fresh eyes. 
``Our course objective, and the real beauty of mathematics,'' I tell them, ``is using specific instances of a problem to understand how it works in full generality.''
And that is exactly what we did.

As mathematicians, we often look to “extreme cases” of a problem to gain intuition for a general solution. 
Mathematics education within the prison system is perhaps the most extreme case, and I believe it sheds light on how we can teach more effectively in a university classroom. 

\section{Satisfying Constraints}

In addition to the usual concerns about teaching a math course, teaching math in prison comes with a unique set of challenges. 
With respect to the aforementioned issue of bringing class materials inside, there are certain canonical classroom resources that are prohibited or impossible to obtain. 
There is no computer, internet access, or even a desk for lecture materials, just a small white board at the front of the room that, just like in university classrooms, is unreliably equipped with dry erase markers. 
There are no office hours, recitations, or means of communicating with students outside of class time. 
If students miss class for a variety of legitimate reasons, like having to work late at their jobs in prison or not being allowed out of their dormitories, they have no way to access material or get in touch with me until the next class. 

During the summer, for example, I taught an elective course called ``The Art of Mathematics,'' in which students investigated several math topics, like algorithms and infinite set cardinalities, that served as inspiration for their own art work. 
TPEI provided sketchbooks, folders, and colored pencils, but I was not permitted to leave the colored pencils with the students when class was over. 
The other class materials required certain stickers to denote that TPEI had given them to these students, in hopes they would not accidentally be confiscated. 

I also brought in art books and textbooks that I had at home to pass around the room, which I was told to collect at the end of each class.
Though I did not, I was very tempted to ignore this rule on two occasions. 
Once when an older student was poring through \textit{The Math Book} by Clifford A. Pickover and asked me excitedly if she could check it out like a library book. 
Then again when a younger student asked me if she could borrow Baby Rudin because she was curious about analysis proofs; she'd also asked me earlier that evening which prerequisite textbooks or classes she'd need in order to teach herself analysis some day.
Later in the course, I learned the latter student was actually studying astronomy in college before her arrest; unfortunately, she'd been ousted from the major due to her grades in the required math classes. 

The culmination of the math art elective course was to be a gallery night, a classroom exhibition of their work on the last day. 
Together, we would reflect on the various topics we'd covered, discussing the mathematical ideas and artistic choices that spoke to us in each other's work. 
I'd invited the TPEI program coordinators, Max and Chloe, who were also excited to see what the students had created.

When we arrived at the classroom, we waited for an unusually long time before a couple of students came in and mentioned they'd had a difficult time getting to our room.
Around 6:45, when Max went back to the central guard station to ask if all TPEI students had been called for class, we were told that a fight had broken out in one of the dorms, and it was now on lockdown.
As a result, despite no involvement in the fight, about half of my students were unable to come and present the art pieces they'd worked on all summer, and by the time the ones who could make it finally arrived, we only had about 45 minutes left.  
Everyone involved was disappointed, but those of us in attendance tried to make the best out of the remainder of the evening.
The exhibit was breathtaking, and the influences and creativity in every piece were truly awe-inspiring.
Max and Chloe were blown away and agreed that the students' exhibit looked and felt professionally curated. 

\begin{figure}[t!]
    \centering
    \begin{subfigure}[t]{0.24\textwidth}
        \centering
        \includegraphics[height=1.85in]{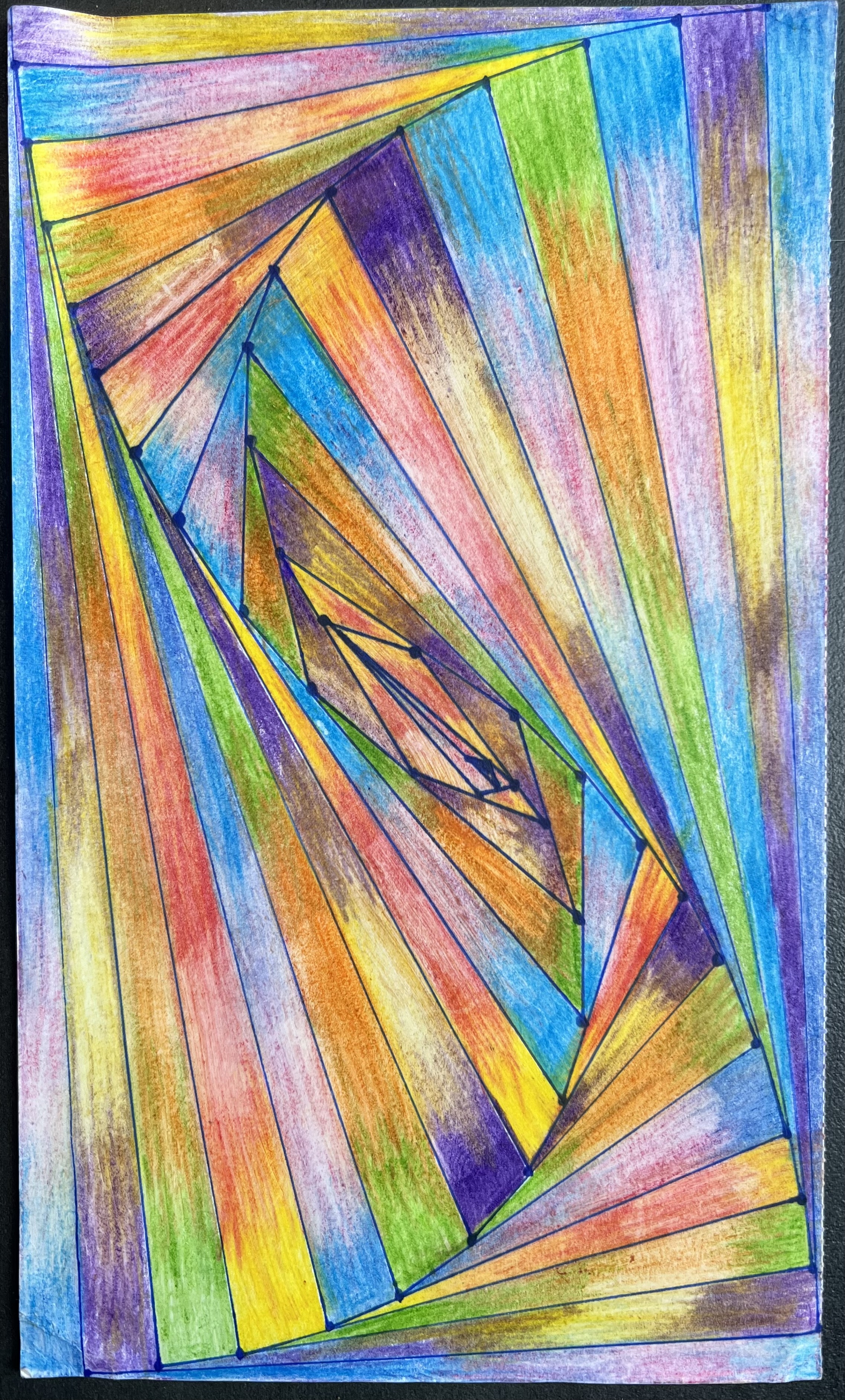}
        \caption{By Tricia B.}
        \end{subfigure}
    \begin{subfigure}[t]{0.24\textwidth}
        \centering
        \includegraphics[height=1.85in]{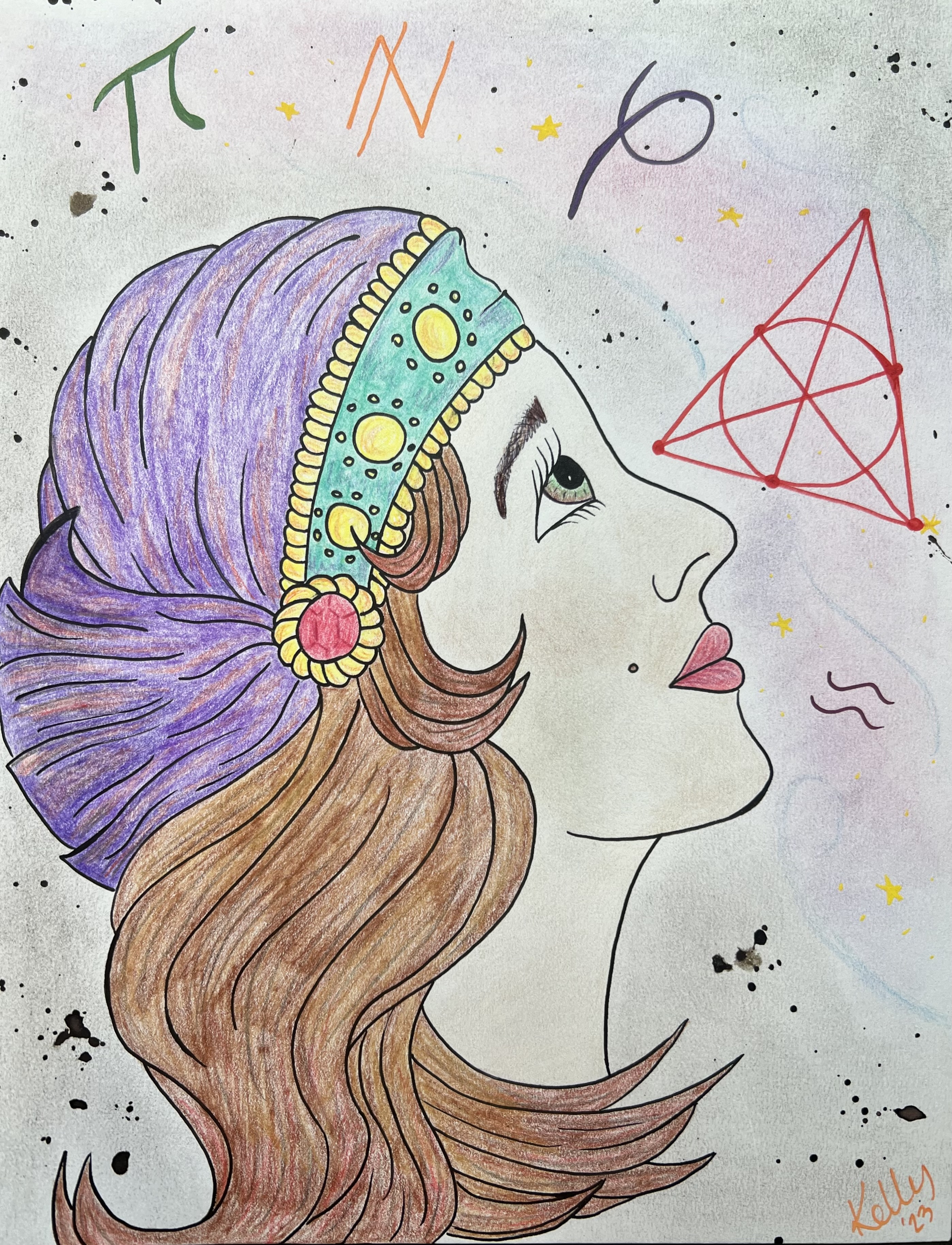}
        \caption{By Kat K.}
    \end{subfigure} \\
            \begin{subfigure}[t]{0.24\textwidth}
        \centering
        \includegraphics[width=1.35in]{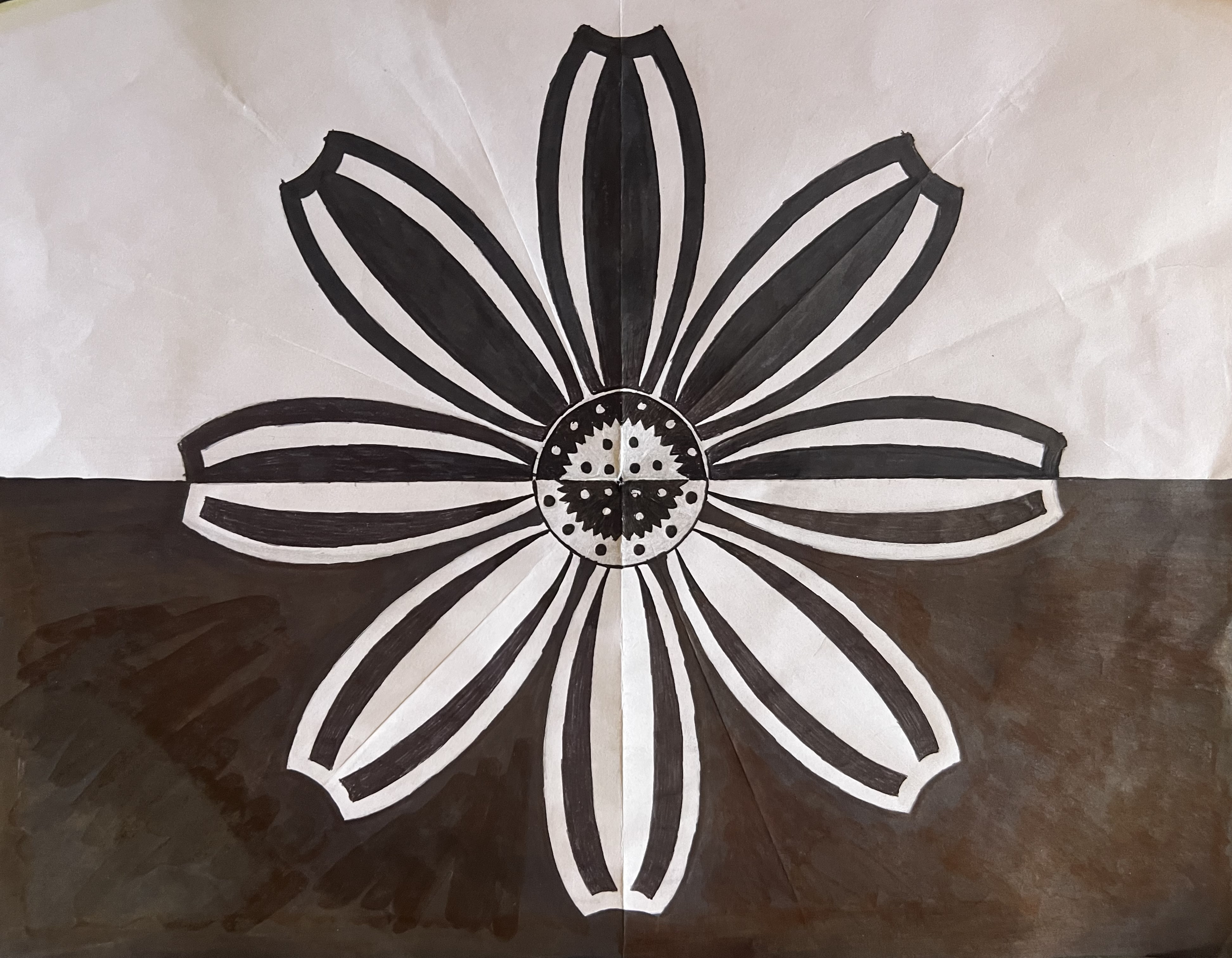}
        \caption{By Marlena R.}
    \end{subfigure}%
    \begin{subfigure}[t]{0.24\textwidth}
        \centering
        \includegraphics[width=1.25in]{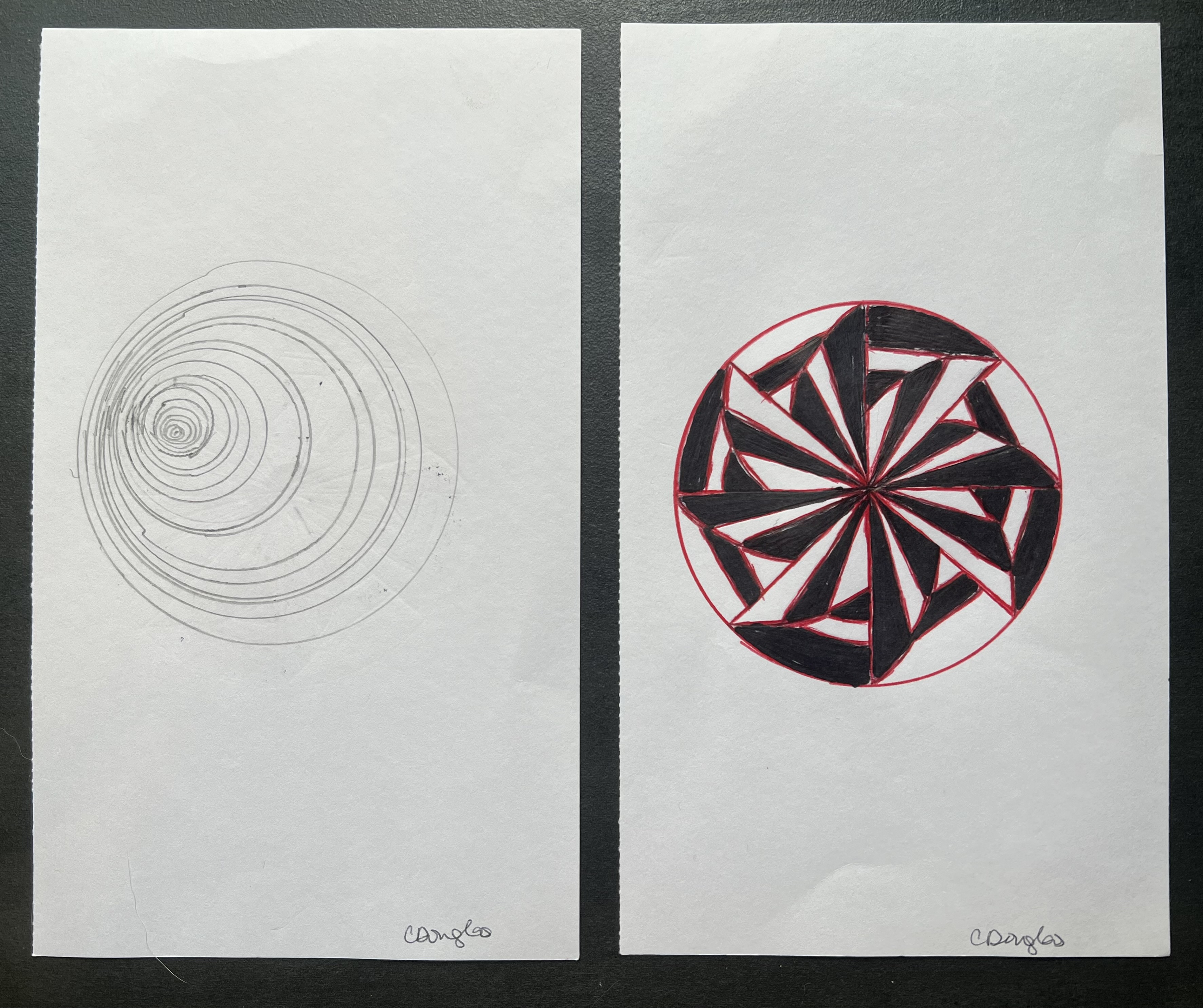}
        \caption{By Tonya D.}
    \end{subfigure} \\
      \begin{subfigure}[t]{0.24\textwidth}
        \centering
        \includegraphics[height=2in]{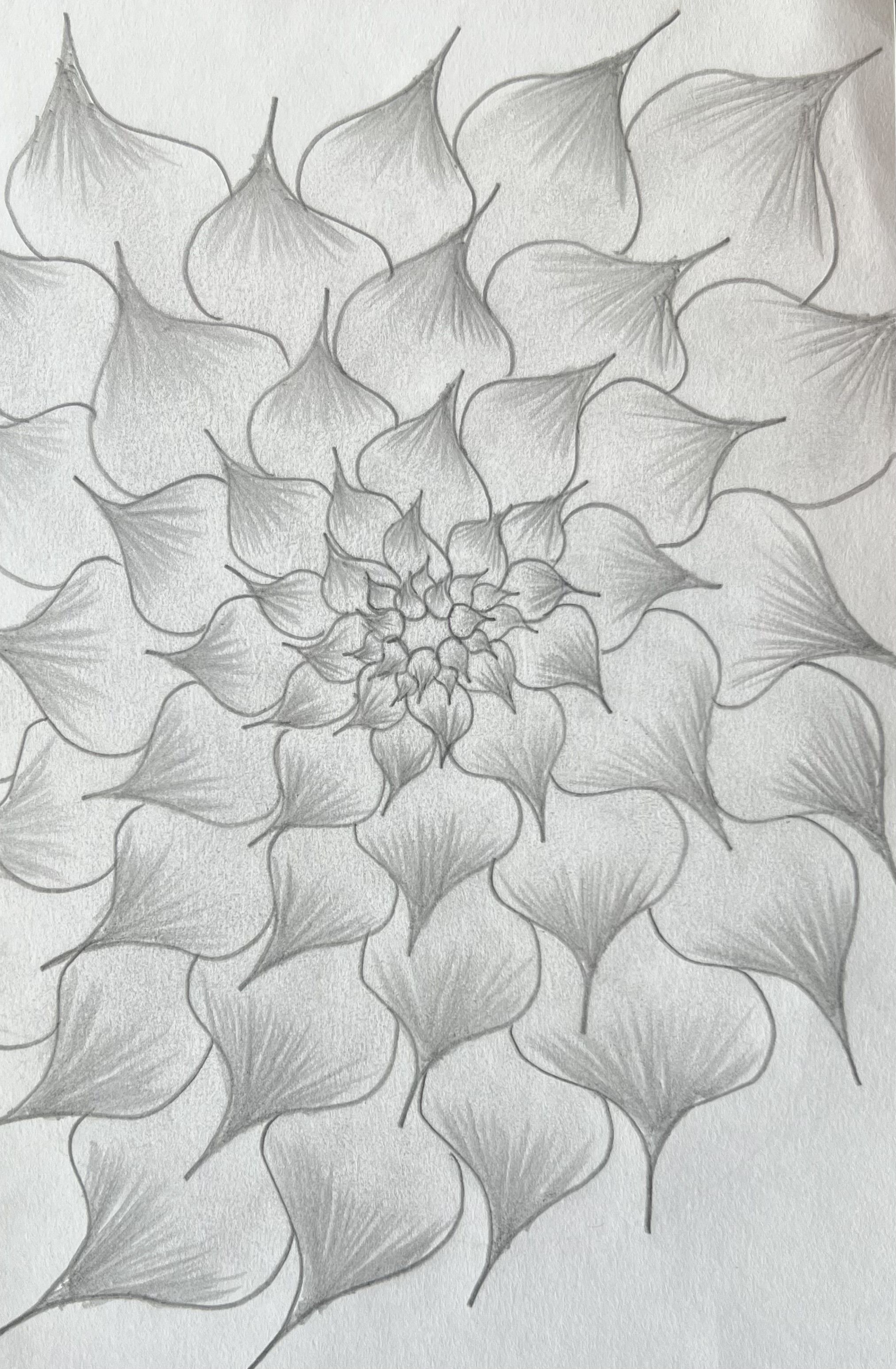}
        \caption{By Tiffany G.}
    \end{subfigure}%
    \begin{subfigure}[t]{0.24\textwidth}
        \centering
        \includegraphics[width=1.4in]{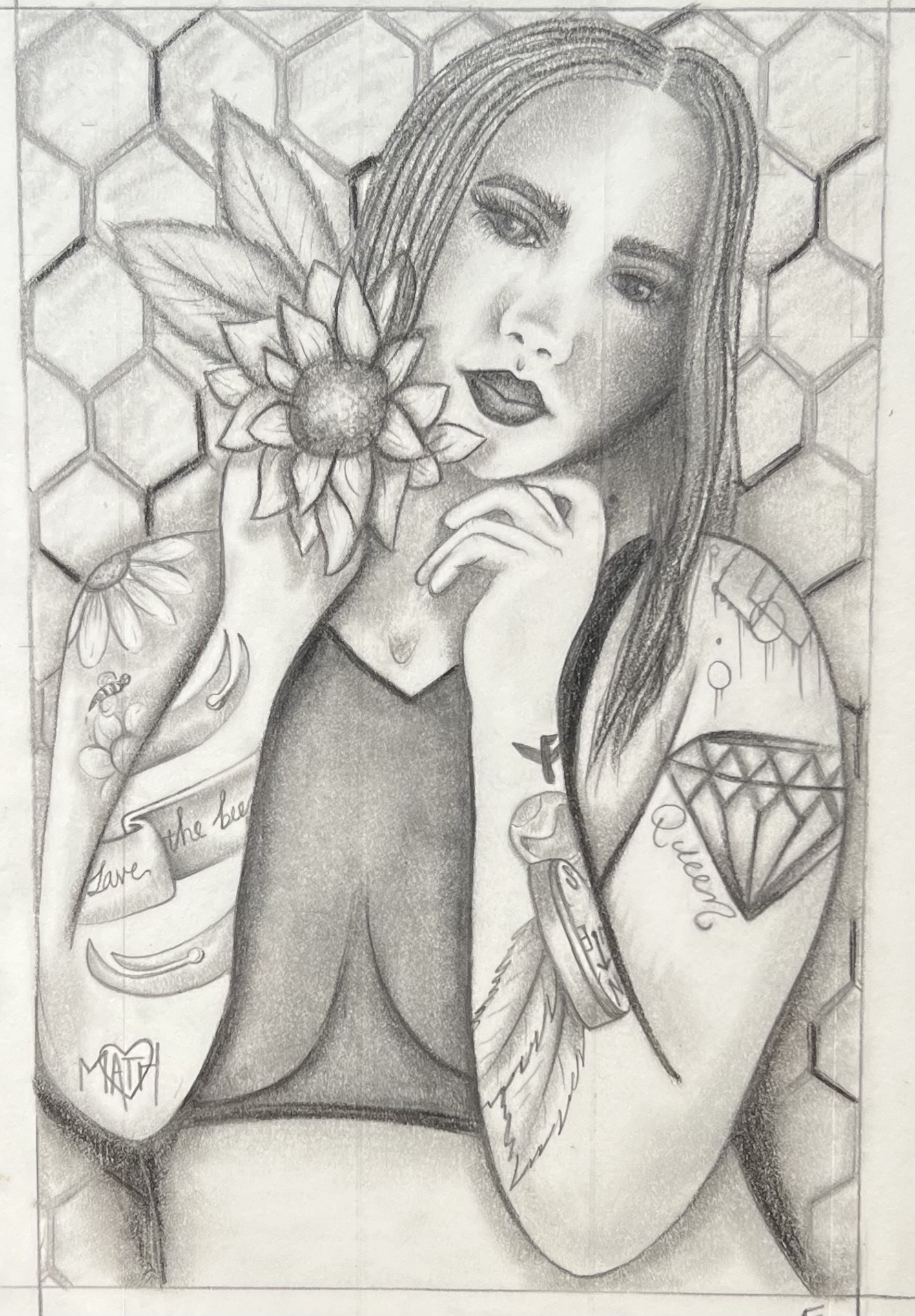}
        \caption{By Amber F.}
    \end{subfigure} \\
    \begin{subfigure}[t]{0.24\textwidth}
        \centering
        \includegraphics[width=1.4in]{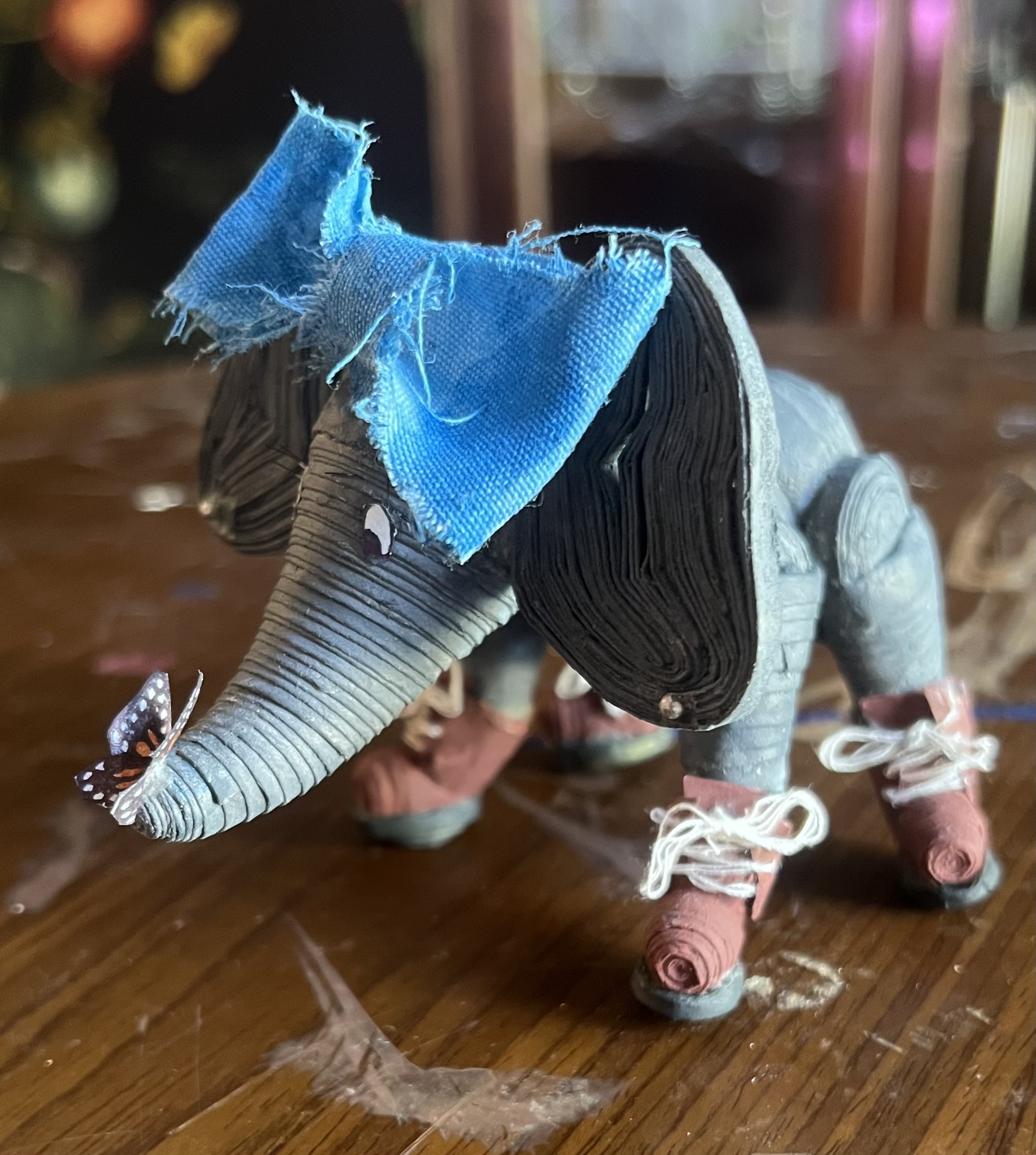}
        \caption{By Heather H.}
    \end{subfigure}%
    \begin{subfigure}[t]{0.24\textwidth}
        \centering
        \includegraphics[width=1.1in]{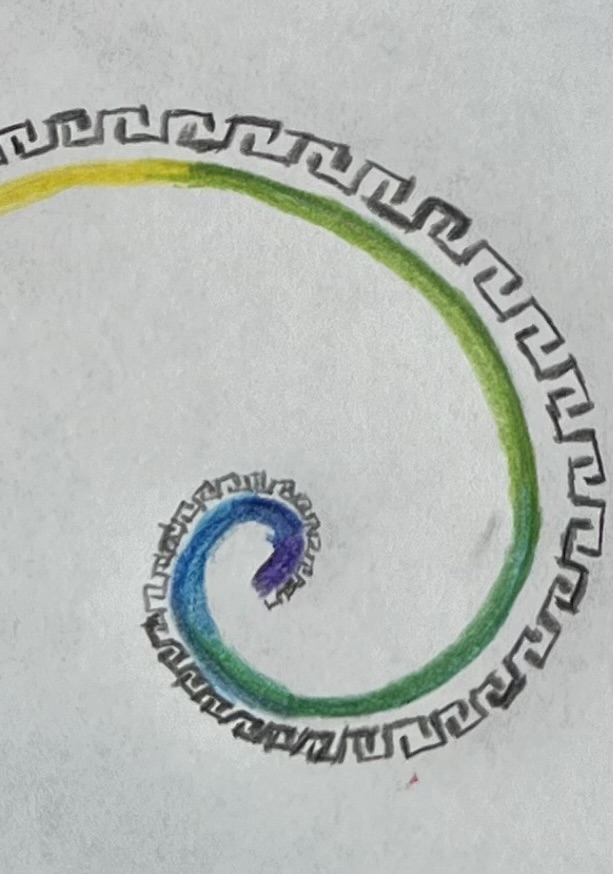}
        \caption{By Erica N.}
    \end{subfigure}%
    \caption{Students' math-inspired art work.}
    \label{fig:art}
\end{figure}

Since we couldn't bring cameras inside to take pictures, my proposal to the students was, with their permission, to borrow their sketchbooks overnight, take photos of their work, then return their sketchbooks to them the next evening at another TPEI instructor's class.  
That way, students wouldn't have to tear out pages from their sketchbooks and would only be without them for less than a day.
However, I was honest about my worry that there could be some unforeseen issue getting the sketchbooks back to them, which was always a risk, and I told them I understood completely if they wanted to hold on to their work just in case.
I was moved almost to tears when every student there that night actually gifted me some of their art work, tearing pages out of their sketchbooks for me to take home and keep.
Those pieces are now framed and displayed in my school office, and Figure~\ref{fig:art} contains photos of their art work. 
In (a), we see a sketch, created with a straight-edge, colored pencils and deoderant for shading, inspired by M.C. Escher's work with H.S.M. Coxeter on the ``limit of infinite smallness'' \cite{Wieting}. 
Figure~\ref{fig:art}(b) shows a self-portrait of the aspiring astronomy student contemplating math topics from the course. 
In (c), the student illustrates the concept of yin and yang using precise mathematical symmetry. 
In Figure~\ref{fig:art}(d), the sketch on the left is a play on the Poincar\'e disk, using circular objects of different radii that the student collected from around the prison. 
On the right in Figure~\ref{fig:art}(d), she depicts her interpretation of Escher's limit of infinite smallness. 
In (e), we see an interpretation of Fibonacci numbers arising in nature, and in (f), a self-portrait of the student in front of a tiled background, with a new (imaginary) tattoo showcasing her love of math on her right elbow. Finally, the photo in (g) showcases a 3-D quilling sculpture of an elephant to represent me, named ``Lil Kate,'' made completely out of small paper strips, dyed with food coloring, and adhered together with a mixture of coffee creamer and water, and in (h), a proportionally-accurate illustration of the Fibonacci spiral.

\section{Developing Intuition}

One of the biggest challenges of teaching in prison is being diligent in my language choices.
I actively try to avoid using idioms or explanations that could trigger mental blocks for students, especially about math subjects in which they'd experienced difficulty or unsuccessful learning outcomes.  
Since I am hoping to build up their intuition on a new mathematical foundation, it's important not to repeat the same explanations that caused them confusion the first time.
Moreover, because of the already difficult environment we are in, I don't want to say something that would cause students any more distress.

A small but pervasive example of being more conscious of my language in the classroom was swapping out the term ``homework'' for ``assignment'' in my college algebra course. 
But I also wanted to make sure my writing achieved the same objectives as my classroom language and tone. 
So every week, I wrote lecture notes that sounded exactly like I would speak during class, while still including the usual textbook definitions, examples, formulas, etc.---and 
 while trying to avoid anything that would elicit a mental block toward the subject.
This style of lecture notes was especially important in situations where students missed class. 
Since there are no dedicated times like office hours, appointments, or recitations for them to get help, I tried to make my notes reminiscent of our classroom: thorough, but not dry, while still containing the explanations they need to understand the material.

Out loud and in text, before introducing any new topic in the course, I wanted to provide motivation for why we were talking about it, but really address why they should care about it. 
As a stereotypical example, I did not have to ask how my students felt about fractions; they all groaned loudly when I said, ``Today, we'll be talking about fractions.'' 
I had anticipated this, mostly because of certain family members and friends who have had an identical response to dealing with fractions. 
I jokingly threw up my hands and acquiesced over the groans: ``Okay, fine, we'll talk about one of my favorite things in math first instead \textellipsis prime numbers!''
We then detoured through prime numbers, prime factorizations, division trees, greatest common divisors, and least common multiples before we finally circled back to fractions. 
Suddenly, the two biggest obstacles to them mastering fractions previously (adding and simplifying) were reduced to problems they had just tackled in a very different setting. 
Once we'd gotten around the initial hurdle to approaching a subject they'd already convinced themselves they'd ``never understand,'' my students no longer experienced the same mental block toward it. 
Instead, I witnessed them appreciate the power of abstraction, to the point that one of them exclaimed, ``Holy crap, I can actually help my kids with their homework now!'' after adding three fractions with different denominators using their least common multiple.

By the time we got to solving linear equations and word problems, the students had a completely different outlook on the material.
After hearing it so much on the first day, I had banned the phrase ``bad at math'' from our classroom, but even if I hadn't, that was not how the students felt anymore.
We had spent multiple class periods developing intuition for the properties of addition, subtraction, multiplication, and division of real numbers, from an abstract perspective but also with the concrete example of debt to understand computations with negative numbers. 
What I'd worried would seem boring to the students was actually what initially piqued their curiosity: we took a specific instance of a problem, then abstracted it away with algebraic tools to study the problem in full generality.
Once they'd seen a given topic from these two different angles, they also had plenty of ``low-stakes'' opportunities to practice.
I frequently reminded them that it's normal to make mistakes when trying something new.
So by the time I said to them, ``Today, we are going to solve some word problems,'' the students were no longer groaning, or feeling anxious about using variables to represent unknowns; their confidence in their abstract problem-solving abilities had been strengthened by the time we spent in the low-stakes material, building a solid foundation from which to work.
Upon solving one of the linear systems arising from a word problem, one of my students proudly announced, ``I feel like a mathematician!''
``You are,'' I replied. 

When I taught the art elective course, I ascribed to a similar ideology.
Each mathematical topic (algorithms, proportion, infinity, and abstraction) was paired with complementary art work from different time periods and regions of the world to illustrate the idea. 
It was surprising how much of the material, which I'd created with the intention of introducing them to advanced undergraduate math they wouldn't have seen before, illuminated other mathematical concepts for them that I hadn't even anticipated. 
When we talked about proportions, for instance, we discussed how the ancient Greeks thought the ideal relationship between the width $w$ and height $h$ of a building is given by the ``Golden proportion''
\begin{align*}
    \varphi = \frac{w}{h} = \frac{w+h}{w}.
\end{align*} 
After we talked about its relationship to the Fibonacci sequence, I mentioned offhandedly that nowadays we can solve explicitly for the Golden ratio $\varphi$ using the quadratic formula; a student immediately raised her hand and asked to see how that would work.
After I showed them the trick of setting $h=1$, they were amazed to see an ``elementary'' formula they'd rotely memorized in school being used to answer such a seemingly unrelated question. 
Later in the same class, when we were discussing the Archimedean method of approximating $\pi$, I wrote the more familiar equation $C = \pi d$ as $\pi = C/d$ to emphasize the ratio.
To my surprise and dismay, the vast majority of them had never seen or considered $\pi$ as a ratio, and it was emotional to watch them eagerly discovering that amazing property of circles every mathematician discovers at some point in their career. 
I began to think deeply about how we are motivating these concepts for students when they first encounter them.

My students' comments in class cast many mathematical concepts in a new light for me as well. 
For example, in talking about algorithmic constructions with \textit{girih} tiles \cite{LuStein07}, we discussed how it would be possible with rudimentary tools to fabricate these tiles with precise angular measurements and arrange them into astonishingly intricate designs.
We talked about the idea of decomposing the polygonal tiles into triangles, and how an understanding of triangles unlocks a lot of possibilities for the tiles' fabrication and design. 
One of the students mused casually, ``So that's why they teach an entire class about triangles. They're like the atoms of shapes.''
I've since begun borrowing that phrase to motivate the subject of trigonometry.

During another class period, we talked about the cardinality of the natural numbers.
I chuckled to myself as I listened to a familiar debate among the students. ``Why would 0 be included? You don't count anything with zero fingers.'' 
Fortunately, they were satisfied that it didn't much matter with respect to cardinality after we wrote out the bijection between the natural numbers and the integers. 
However, that awareness came back to bite me in the form of a deeply profound and unexpected question that I received while explaining Cantor's diagonalization argument.
I had just demonstrated his proof by contradiction: if we try to enumerate all of the elements $s_1, s_2, s_3, \ldots$ of the set $T$ of infinite sequences of 0's and 1's, it is always possible to construct an element $s \in T$ that differs from $s_k$ in the $k$th position, so that element of $T$ actually wasn't enumerated in our list. 
``But why can't you just call that element  $s_0$? Then wouldn't you be able to count them all with counting numbers since it doesn't matter if we include 0?''
The question was so subtle and clever that it caught me off guard.
``You're thinking like a mathematician,'' I replied, before we spent the next few minutes verifying that Cantor had indeed gotten it right, though probably not on his first try.

\section{Solving Problems}

Once I was talking about my experience teaching math in prison with a friend of mine who was formerly incarcerated. 
He'd served a ten-year sentence over the entirety of his 20s, during which time he had the opportunity to take several math courses for high school degree equivalency.
He's now earning his bachelors degree while working full-time as an operations manager at a water treatment facility.
While math courses equipped him with necessary skills for his new career, he credits those classes with something even more important.
As a creative writer, my friend had always felt more passionate about writing classes in school.
He admitted to me that he'd never liked that there was ``only one right answer'' in his math classes.
But once he was incarcerated, solving math problems became a mentally stimulating and comforting activity.
With so much time to think and reflect, he realized how many problems in life lack a clear-cut solution, and he began to appreciate the existence and uniqueness of the solutions to the problems in his math assignments.

I told my students about this conversation with my friend, and I asked them if they felt similarly about learning math while incarcerated.
Several of them  participate in entrepreneurial and vocational training programs, and those students echoed the importance of math for their new career paths.
In fact, one student  made parole during our spring algebra class, and she still wanted to finish the course remotely post-incarceration to help  her earn her EMT certifications.
(This situation is one of many in which the TPEI administrative team of volunteers is essential to the program's success.) 

Several other students said that my friend's point about taking comfort in having a right answer really resonated with them, even though they'd had terrible previous experiences with math. One of them was the student who'd had to drop her astronomy major because of her math classes, who'd asked about teaching herself Baby Rudin. 
Another student had been concurrently earning her high school diploma during our class (``only 50 years late!'' as she'd often say), and she explained to me how learning math again had helped her find balance in her life. 
``Everything in math has an opposite,'' she said, ``I am a `yin and yang' type of person, and I guess that's why I love math. You can always add back anything you subtract, or multiply anything you divide.''
Her friend chimed in, ``Plus, even just knowing there are infinitely many possible choices to start from gives you options for all of the gray areas in life.''
Yet another student, who was nervous about her upcoming parole hearing the following week, reiterated that being able to solve a problem and find a correct answer out of infinitely many possible choices was going to be vital to her success post-incarceration.
Like so many of her  classmates, she'd also been interested in STEM growing up; she wanted to be an astronaut when she entered high school, before dropping out because of math. 
``I'm nervous about finding a job and a permanent place to live, but I have people to help me for now,'' she said.
``I probably can't be an astronaut anymore because of my felony charge, but I still want to enroll in college and earn my B.S. once I'm back on my feet.''
In spite of all the problems and uncertainties that she faced post-incarceration, she'd already begun identifying solutions.

On the last day of the art elective class, after the exhibition, I announced that I would be teaching a brand new course in the fall that had never been offered by TPEI. 
It is a credit-bearing math course through UT Austin that is a prerequisite for UT's calculus sequence. 
My student---the aspiring astronomer who'd been enthralled by Baby Rudin---came up to talk to me after the other students had said their thank-you's and goodbye's. 
``I really want to take the course, but I've failed this subject before and I'm worried I will again,'' she said nervously.
``I haven't done math in a long time.''
Thinking back over all of the insights she'd shared during the course, I smiled and reassured her that I knew she could do it and that I would be there to help her.
``Besides,'' I replied, ``you've been thinking like a mathematician this whole time.'' 

\bibliography{main.bbl}

\begin{bibdiv}
\begin{biblist}

\bib{LuStein07}{article}{
      author={Lu, Peter~J.},
      author={Steinhardt, Paul~J.},
       title={Decagonal and quasi-crystalline tilings in medieval islamic
  architecture},
        date={2007},
     journal={Science},
      volume={315},
       pages={1106\ndash 1110},
}

\bib{Wieting}{article}{
      author={Wieting, Thomas},
       title={Capturing infinity: The circle limit series of m.c. escher},
        date={2010March},
     journal={Reed Magazine},
}

\bib{PPI22}{article}{
      author={Wang, Leah},
      author={Sawyer, Wendy},
      author={Herring, Tiana},
      author={Widra, Emily},
       title={Beyond the count: A deep dive into state prison populations},
        date={2022},
      eprint={www.prisonpolicy.org/reports/beyondthecount.html},
}

\end{biblist}
\end{bibdiv}

\end{document}